\documentclass{article}
\usepackage{amssymb}
\usepackage{amsmath}
\usepackage{latexsym}
\usepackage{mathrsfs}%花体
\usepackage{graphicx}%插图
\usepackage[all]{xy}%交换图
\newtheorem{Th}{Theorem}
\newtheorem{Lemma}[Th]{Lemma}
\newtheorem{Coro}[Th]{Corollary}

\title{Schur Polynomials through Lindstr\"om Gessel Viennot Lemma}
\author{XIONG Rui}
\begin{document}

\maketitle

\begin{abstract}
  In this article, we use Lindstr\"om Gessel Viennot Lemma to give a short, combinatorial, visualizable proof of the identity of Schur polynomials
  --- the sum of monomials of Young tableaux equals to the quotient of determinants. As a by-product, we have a proof of Vandermonde determinant without words. We also prove the cauchy identity. In the remarks, we discuss factorial Schur polynomials, dual Cauchy identity and the relation beteen Newton interpolation formula.
\end{abstract}

%\tableofcontents
\section{Introduction}

The main purpose of this article is to give a combinatorial proof of the identity of Schur polynomials
$$\sum_{T \textrm{ Young tableaux}} x^T=
%\frac{\det\left(\begin{matrix}
%x_1^{\lambda_1+n-1} & x_2^{\lambda_1+n-1} & \cdots & x_n^{\lambda_1+n-1}\\
%x_1^{\lambda_2+n-2} & x_2^{\lambda_2+n-2} & \cdots & x_n^{\lambda_2+n-2}\\
%\vdots & \vdots & \ddots & \vdots \\
%x_1^{\lambda_n} & x_2^{\lambda_n} & \cdots & x_n^{\lambda_n}\\
%\end{matrix}\right)}{\det\left(\begin{matrix}
%x_1^{n-1} & x_2^{n-1} & \cdots & x_n^{n-1}\\
%x_1^{n-2} & x_2^{n-2} & \cdots & x_n^{n-2}\\
%\vdots & \vdots & \ddots & \vdots \\
%x_1^{0} & x_2^{0} & \cdots & x_n^{0}\\
%\end{matrix}\right)}
\frac{\det(x_i^{\lambda_j+n-j})}{\det(x_i^{n-j})} $$
using this Lindstr\"om Gessel Viennot lemma. The classic algebraic proof,
see Fulton and Harris \cite{zbMATH00051906} page 462, proving $\det(h_{\lambda_i+i-j})=\frac{\det(x_i^{\lambda_j+n-j})}{\det(x_i^{n-j})}$.
The classic combinatorial proof (analysing a bijection carefully), see Mendes, and Remmel \cite{mendes2015counting}, page 40.

\paragraph{Lindstr\"om Gessel Viennot Lemma. }
Let $\Gamma=(V,E)$ be a locally finite directed acyclic graph.
Let $\{X_e:e\in E\}$ be a set of pairwise commutative indeterminants.
For each path $P:v_0\stackrel{e_1}\to v_1\to \cdots \stackrel{e_k}\to v_k$, we define $X_P=X_{e_1}\cdots X_{e_k}$.
For two vertices $u,v\in V$, we define $e(u,v)$ to be the sum of $X_P$ with $P$ going through all path from $u$ to $v$.

Now, fix some integer $n>0$, and two subsets of order $n$ of $V$, say $A=\{a_1,\ldots,a_n\}$ and $B=\{b_1,\ldots,b_n\}$.
We have the following lemma.

\begin{Lemma}[Lindstr\"om Gessel Viennot Lemma, \cite{gessel1989determinants}]\label{LGVLemma} As notations above,
$$%\det\left(\begin{matrix}
%e(a_1,b_1) & \cdots & e(a_1,b_n)\\
%\vdots & \ddots & \vdots \\
%e(a_n,b_1) & \cdots & e(a_n,b_n)
%\end{matrix}\right)
\det(e(a_i,b_j))=\sum_{\sigma\in \mathfrak{S}_n} \operatorname{sgn}\sigma\sum_{\{P_1,\ldots,P_n\}} X_{P_1}\cdots X_{P_n}$$
where $\{P_1,\ldots,P_n\}$ goes through all pairwise non-intersecting paths with each $P_i$ from $a_i$ to $b_{\sigma(i)}$.
\end{Lemma}
%
%The proof is not hard, consider the following operator on intersecting paths.
%$$\includegraphics{Operator}, $$
%Find the lowest and firstmost place the intersection happens, then this operator will pair and cancell all of them.
%
%The case when there is \emph{only one} possible permutation $\operatorname{id}\in \mathfrak{S}_n$ without intersecting will be mainly used.
%%For example, \cite{Fulmek2010ViewingDA} discussed a lot of applications of this theorem.

\paragraph{Schur polynomials. }
For $n>0$, the $n$-variable \emph{Schur polynomial} of a partition $\lambda: \lambda_1\geq \cdots\geq \lambda_n$ is defined to be
$$S_{\lambda}=\sum_{T} x^T,$$
where $T$ goes through all Young tableaux (weakly increasing in row and strictly increasing in column)
of alphabet $\{1,\cdots,n\}$ with shape $\lambda$. If we denote $N_i=\# T^{-1}\{i\}$, then $x^T$ is defined to be
$x_1^{N_1}\cdots x_n^{N_n}$.

We know the following identity.
\begin{Th}[Jacobi-Trudy identity]\label{JacobiTrudyidentity} As notation above,
$$S_{\lambda}=\det(h_{\lambda_i+i-j}),$$
where $h_k$ is the sum of all monomials of degree $k$ in $n$ variables.
\end{Th}

One of the combinatorial proof is to use Lindstr\"om Gessel Viennot lemma \ref{LGVLemma} above.
Since we will use the proof of it, let me state the proof briefly.

%
%\begin{figure}[h]
%  \centering
%  \includegraphics{Jacobi-Trudy}\\
%  \caption{An example}\label{eg1}
%\end{figure}

Consider the \emph{lattice graph} $\Gamma=(V,E)$, with $V=\mathbb{Z}\times \mathbb{Z}$ the plane lattice points,
and $E$ all the edge from $(i,j)$ to $(i,j+1)$ and $(i+1,j)$. Assign the edge $(i,j)\to (i+1,j)$ by weight $x_j$, the rest by $1$,
and consider
$$\begin{array}{l}
A=\{a_1=(1,1),a_2=(2,1),\ldots,a_n=(n,1)\};\\
B=\{b_1=(1+\lambda_n,n), b_2=(2+\lambda_{n-1},n),\ldots,b_n=(n+\lambda_1,n)\}.
\end{array}$$
Then the left hand side of Lindstr\"om Gessel Viennot Lemma is exactly $\det(h_{\lambda_i+i-j})$,
and right hand side is exactly Schur polynomial. % see figure \ref{eg1} for a clue.
The proof is contained in wikipedia, \cite{wikiLGVLemma}.
For a complete proof, see for example Prasad \cite{Prasad2018AnIT}.  %Mendes, and Remmel \cite{mendes2015counting}, page 61.

\section{Main Lemma}

Here we consider the same graph $\Gamma=(V,E)$, the lattice graph, as above but of different weight.
More precisely, assign the edge $(i,j)\to (i+1,j)$ by $x_j-x_{i+j}$, and the rest by $1$.
%$$\xymatrix@!C=3pc{
%\vdots&\vdots&\vdots&\\
%\bullet\ar[r]|{(x_3-x_4)}\ar[u] & \bullet \ar[r]|{(x_3-x_5)}\ar[u] & \bullet \ar[r]|{(x_3-x_6)}\ar[u]&\cdots\\
%\bullet\ar[r]|{(x_2-x_3)}\ar[u] & \bullet \ar[r]|{(x_2-x_4)}\ar[u] & \bullet \ar[r]|{(x_2-x_5)}\ar[u]&\cdots\\
%\bullet\ar[r]|{(x_1-x_2)}\ar[u] & \bullet \ar[r]|{(x_1-x_3)}\ar[u] & \bullet \ar[r]|{(x_1-x_4)}\ar[u]&\cdots
%}$$

\begin{Lemma}Let $a=(1,1)$ and $b=(m,n)$, with $m,n>0$, then
$$e(a,b)=(x_1-x_{m+n-1})\cdots (x_1-x_{n+2})(x_1-x_{n+1}). $$
We use the convention that when $m+n-1<n+1$ the above expression equals to $1$.
\end{Lemma}
This follows from induction. Firstly, the expression holds when $m=1$ or $n=1$.
Generally, if $m,n>1$, let $b'=(m-1,n)$ and $b''=(m,n-1)$, as below
$$\xymatrix@!C=8pc{
e(a,b')\ar[r]|{(x_m-x_{m+n-1})} & e(a,b)\\
\cdots & e(a,b'')\ar[u]}$$
Then
$$\begin{array}{rll}
e(a,b) & = e(a,b')(x_n-x_{m+n-1})+e(a,b'')\\
& =(x_1-x_{m+n-2})\cdots (x_1-x_{n+2})(x_1-x_{n+1})(x_n-x_{m+n-1})\\
& +(x_1-x_{m+n-2})\cdots (x_1-x_{n+2})(x_1-x_{n+1})(x_1-x_{n}) & \textrm{by induction}\\
& = (x_1-x_{m+n-1})\cdots (x_1-x_{n+2})(x_1-x_{n+1}).
\end{array}$$
The proof is complete.
\bigbreak

\begin{figure}[h]
  \centering
  \includegraphics{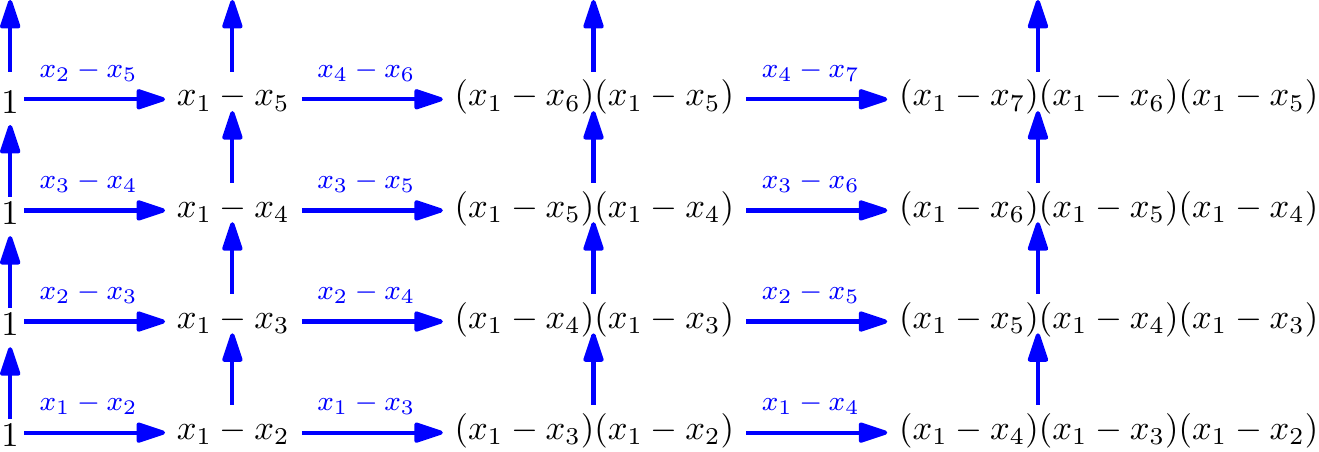}\\
  \caption{The table}\label{mylemmatable}
\end{figure}

See figure \ref{mylemmatable} for first several terms.

\begin{Coro}\label{Corollaryoflemma}
Let $a=(1,t)$ and $b=(m,n)$, with $m>s>0,n>t>0$. Then
$$e(a,b)=(x_t-x_{m+n-1})\cdots (x_t-x_{n+2})(x_t-x_{n+1}). $$
In particular, when $0=x_{n+1}=x_{n+2}=\cdots$, $e(a,b)=x_t^{m-1}$.
\end{Coro}
By replacing $(x_1,x_2,\ldots)$ by $(x_{t},x_{t+1},\ldots)$.

\bigbreak
As an application, we can get the classic Vandermonde determinant from above computation and Lindstr\"om Gessel Viennot Lemma \ref{LGVLemma}.

\begin{figure}[h]
  \centering
  \includegraphics{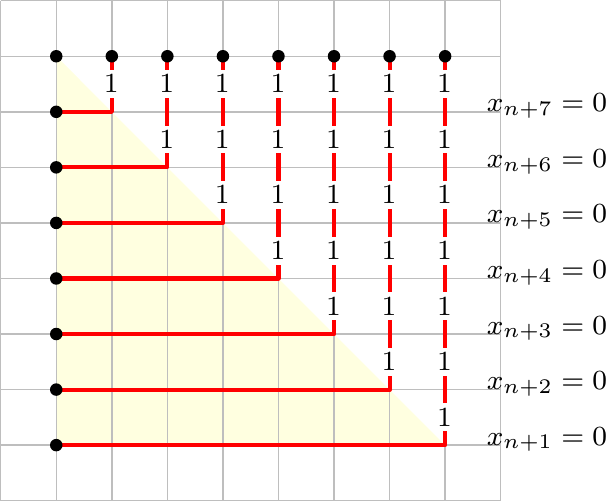}\\
  \caption{Prove Vandermonde determinant without words}\label{Vandermonde}
\end{figure}

\begin{Th}[Vandermonde]\label{Vandermondedet} For $n$ variables, the determinant $\det(x_i^{n-j})=\prod_{i<j}(x_i-x_j)$.
\end{Th}
Consider the graph above, with $0=x_{n+1}=x_{n+2}=\cdots$
and consider
$$\begin{array}{l}
A=\{a_1=(1,1), a_2=(1,2),\cdots, a_n=(1,n)\}\\
B=\{b_1=(n,n), b_2=(n-1,n), \cdots, b_n=(1,n)\}
\end{array}$$
Then by \ref{Corollaryoflemma} above $e(a_i,b_j)=x_i^{n-j}$.
But the only non-intersecting paths are as figure \ref{Vandermonde}.
Apply Lindstr\"om Gessel Viennot Lemma \ref{LGVLemma}, we see $\det(x_i^{n-j})=\prod_{i<j}(x_i-x_j)$.

\section{Schur polynomial by quotient of determinants}

\begin{Th}The $n$-variable Schur polynomial of a partition $\lambda: \lambda_1\geq \cdots\geq \lambda_n$
$$S_{\lambda}=\frac{\det(x_i^{\lambda_j+n-j})}{\det(x_i^{n-j})}. $$
\end{Th}

\begin{figure}[h]
  \centering
  \includegraphics{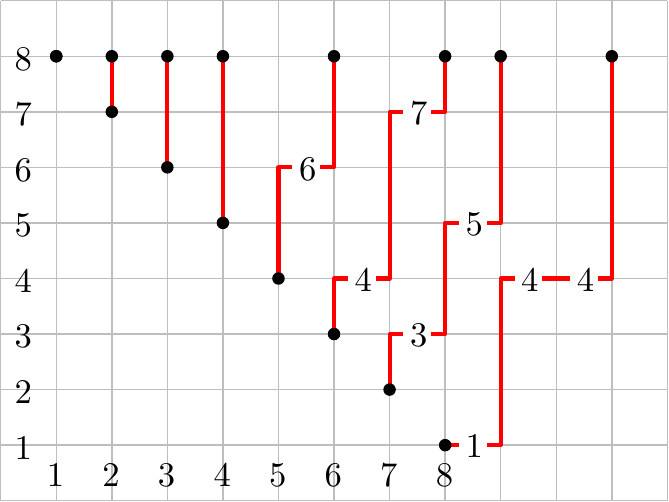}\\
  \caption{First step of reduction}\label{Reduction}
\end{figure}

\textbf{Step 1. } Remind the proof of Jacobi Trudy identity \ref{JacobiTrudyidentity}.
Since the path from $(i,1)$ in any set of non-intersecting paths only went vertically before height $n-i+1$,
so it is equivalent to change $(i,1)$ by $a_i'=(i,n-i+1)$. That is,
$$\det(e(a'_i,b_j))=S_{\lambda}.\eqno{(1)}$$
See figure \ref{Reduction}.

Now, if we apply $0=x_{n+1}=x_{n+2}=\cdots$ in the weight of last section,
then in the region $\{(i,j):i+j\geq n\}$, it coincides with that in the proof of Jacobi Trudy identity \ref{JacobiTrudyidentity}.

\begin{figure}[h]
  \centering
  \includegraphics{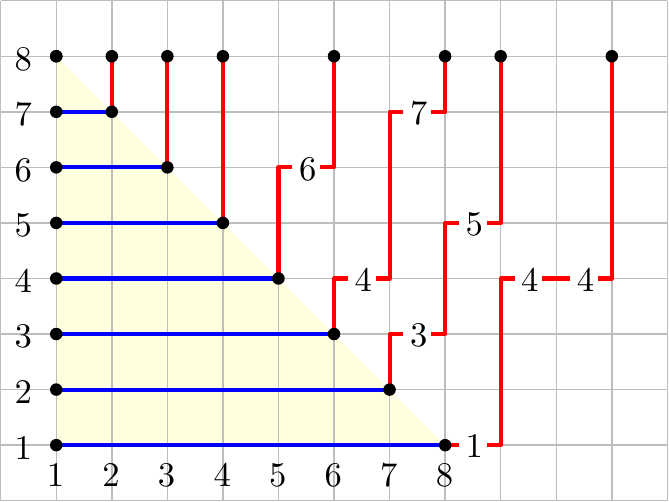}\\
  \caption{Second step of reduction}\label{Reduction2}
\end{figure}

\textbf{Step 2. } Let us consider $a_i''=(1,n-i+1)$. It is easy to see
$$e(a_i'',b_j)=\sum_k e(a_i'',a_k')\cdot e(a_k',b_j), $$
so
$$\det(e(a_i'',b_j))=\det(e(a_i'',a_j'))\cdot \det(e(a_i',b_j)).\eqno{(2)}$$
And trivially,
$$\det(e(a_i'',a_j'))=\prod_{i<j}(x_i-x_j)=\det(x_i^{n-j}), \eqno{(3)}$$
the Vandermonde determinant \ref{Vandermondedet}.
See figure \ref{Reduction2}.

\begin{figure}[h]
  \centering
  \includegraphics{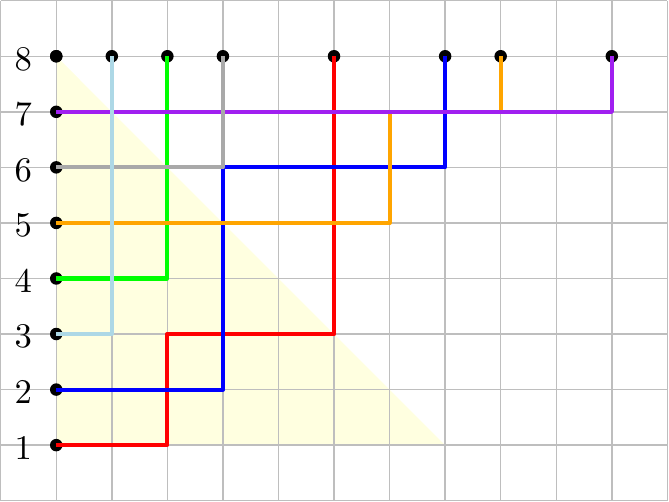}\\
  \caption{Last step of reduction}\label{Reduction3}
\end{figure}

\textbf{Step 3. } By corollary \ref{Corollaryoflemma}, $e(a_{n+1-i}'',b_{n+1-j})=x_i^{\lambda_j+n-j}$, so
$$\det(e(a_i'',b_j))=\det(x_i^{\lambda_j+n-j})\eqno{(4)}$$
See figure \ref{Reduction3}.

Hence the equation $(1,2,3,4)$ proves the assertion.

\section{Cauchy identity}

We can use the identical way to prove Cauchy identity.

\begin{figure}[h]
  \centering
  \includegraphics{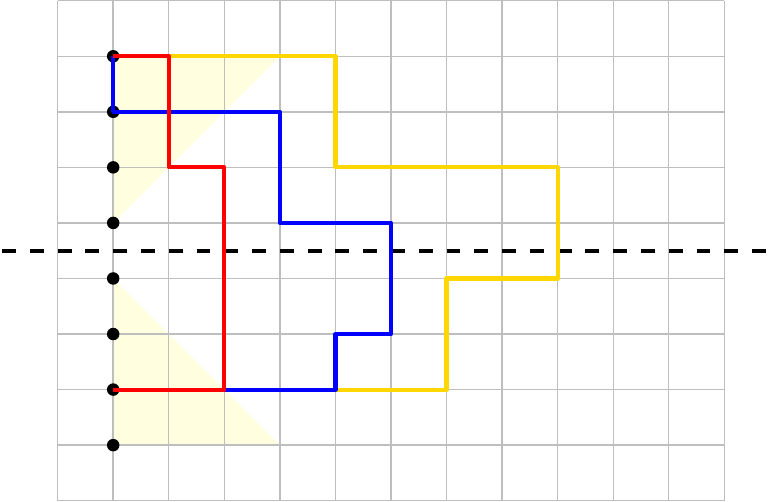}\\
  \caption{Cutting the line $y=n+1/2$}\label{Cuttingtheline}
\end{figure}

\begin{figure}[h]
  \centering
  \includegraphics{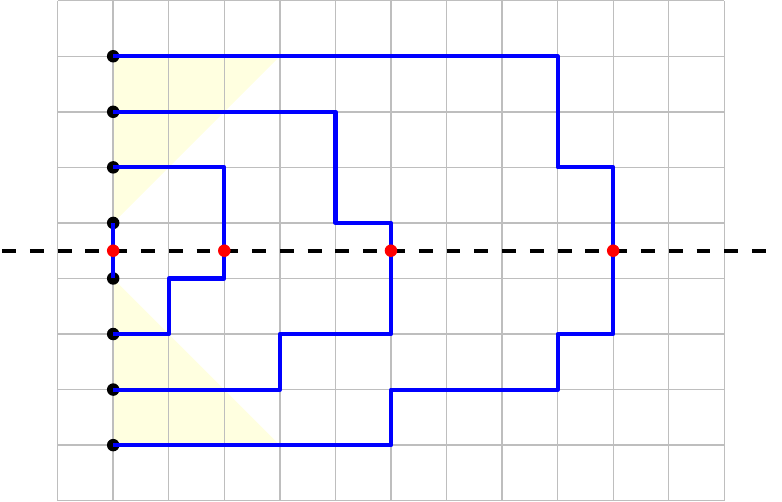}\\
  \caption{Cutting the line $y=n+1/2$}\label{Cuttingtheline2}
\end{figure}

\begin{Th}[Cauchy]\label{CauchyID}For fixed $n>0$, and $x_1,\ldots,x_n,y_1,\ldots,y_n$ two sets of variables,
$$\frac{1}{\prod_{i<j}(x_i-x_j)\prod_{i<j}(y_i-y_j)}\det\left(\frac{1}{1-x_iy_j}\right)=\sum_{\lambda} S_{\lambda}(x)S_{\lambda}(y), $$
with $\lambda$ going through all young diagrams with no more than $n$ rows.
\end{Th}
Here we consider a variant of our graph. Consider the graph $\Gamma=(V,E)$,
$$V=\{1,2,3,\ldots\}\times \{1,2,\cdots,2n-1,2n\}$$
with edge $(i,j)\to (i,j+1)$ with weight $1$,
and $(i,j)\to (i+1,j)$ of weight $x_j-x_{i+j}$ if $j\leq n$,
$(i+1,j)\to (i,j)$ of weight $y_{2n-j}-y_{i+2n-j}$ if $j\geq n+1$.
Then apply $0=x_{n+1}=x_{n+2}=\cdots$ and $0=y_{n+1}=y_{n+2}=\cdots$. Consider the points
$$a_1=(1,1),a_2=(1,2),\cdots,a_n=(1,n)$$
and
$$b_n=(1,n+1),b_{n-1}=(1,n+2),\cdots,b_1=(1,2n). $$
Then cutting the line $y=n+1/2$, by corollary \ref{Corollaryoflemma} we see (see figure \ref{Cuttingtheline})
$$e(a_i,b_j)=\sum_{k=0}^\infty x_i^ky_j^k=\frac{1}{1-x_iy_j}. $$
Apply Lindstr\"om Gessel Viennot Lemma \ref{LGVLemma}.
By the same argument of cut, we see the right hand side is actually $\sum_{\lambda} S_{\lambda}(x)S_{\lambda}(y)$, see figure \ref{Cuttingtheline2}.
The proof is complete.
\bigbreak

Actually, the left hand side of above identity turns out to be $\prod_{i,j=1}^n\frac{1}{1-x_iy_j}$, known as Cauchy determinant.

For the classic proof, see Richard and Sergey \cite{stanleyfomin1999} theorem 7.12.1, using Robinson Schensted Knuth Algorithm.
\section{Remarks}

\paragraph{Factorial Schur polynomials. }
Firstly,
there are also \emph{factorial Schur polynomials}
$$S_{\lambda}(x\mid a)=\sum_{T}\prod_{(i,j)\in \lambda}(x_{T(i,j)}-a_{T(i,j)+i-j})$$
with $A$ going through all Young tableaux of $\lambda$. By a completely same argument, but $x_{n+i}=a_{i}$,
$$S_{\lambda }(x\mid a)=\frac{\det[(x_{j}\mid a)^{\lambda _{i}+n-i}]}{\det(x_i^{n-j})}$$
with $(y\mid a)^k = (y-a_1)... (y-a_k)$.

%\bigbreak

\begin{figure}[h]
  \centering
  \includegraphics{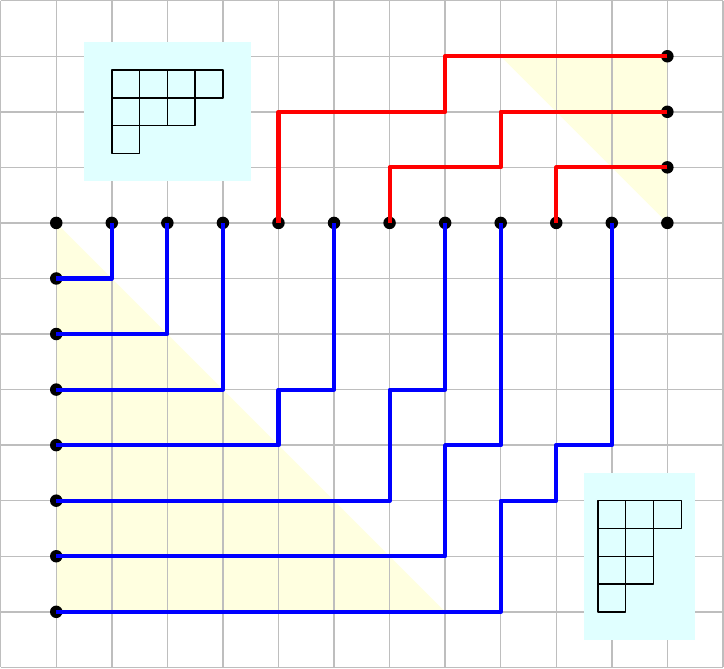}\\
  \caption{Dual Cauchy identity}\label{DualCauchyId}
\end{figure}

\paragraph{Dual Cauchy identity. }
Secondly,
there is also \emph{dual Cauchy identity},
$$\prod_{i=1}^n\prod_{j=1}^m(1+x_iy_j)=\sum_{\lambda} S_{\lambda}(x)S_{\lambda'}(y)$$
where $\lambda$ goes through all young diagrams of at lost $n$ rows and $m$ columns, and $\lambda'$ is the conjugation of $\lambda$.
The classic proof, see Richard and Sergey \cite{stanleyfomin1999} theorem 7.14.3, using dual Robinson Schensted Knuth Algorithm.
This can be proven like what we did in proof of Cauchy identity \ref{CauchyID} but a little careful about the sign. More precisely,
consider the graph $\Gamma=(V,E)$,
$$V=\{1,2,3,\cdots,m+n\}\times \{1,2,3\cdots,m+n-1\}$$
with edge $(i,j)\to (i,j+1)$ with weight $1$,
and $(i,j)\to (i+1,j)$ of weight $x_j-x_{i+j}$ if $j\leq n$,
$(i,j)\to (i+1,j)$ of weight $-(y_{m+n-j}-y_{(m+n-i)+(m+n-j)})$ if $j\geq n+1$. Similar, we apply $0=x_{n+1}=x_{n+2}=\cdots$ and $0=y_{m+1}=y_{m+2}=\cdots$
Then consider
$$(1,1), (1,2),\ldots,(1,n),\quad (m+n,n),(m+n,n+1),\ldots,(m+n,m+n-1)$$
and
$$(1,n), (2,n),\ldots,(m+n,n). $$
See figure \ref{DualCauchyId}.
The sign of permutation cancels the sign of $y_\bullet$'s. Apply Lindstr\"om Gessel Viennot Lemma, \ref{LGVLemma}, the left hand side is
$$\det\left(\begin{matrix}
x_1^{m+n-1} & \cdots & x_n^{m+n-1} & (-y_1)^{0} & \cdots & (-y_m)^0\\
\vdots & \ddots & \vdots & \vdots & \ddots & \vdots \\
x_1^{0} & \cdots & x_n^{0} & (-y_n)^{m+n-1} & \cdots & (-y_m)^{m+n-1}
\end{matrix}\right)\eqno{(*)}$$
and the right hand side is
$$\prod_{1\leq i<j\leq n}(x_i-x_j)\prod_{1\leq i<j\leq m}(y_i-y_j)\sum_{\lambda} S_{\lambda}(x)S_{\lambda'}(y). $$
By a direct computation of $(*)$ (a good exercise of linear algebra), we get dual Cauchy identity.

\begin{figure}[h]
  \centering
  \includegraphics{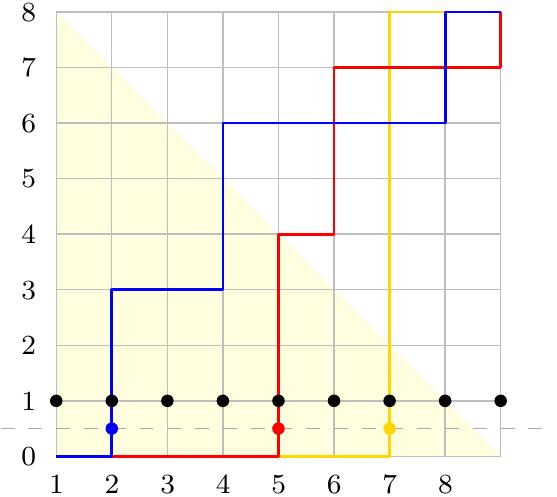}\\
  \caption{Newton interpolation formula}\label{Newtoninterpolation}
\end{figure}

\paragraph{Newton interpolation formula. }
Last but no mean least, our weight chosen is relative to \emph{Newton interpolation formula}. If we view $x_0=x$ as variable,
and $x_1,x_2,\ldots,x_n$ as data points, as usual $0=x_{n+1}=x_{n+2}=\cdots$.
Cutting the path from $(1,0)$ to $(n,n)$ by $y=\frac{1}{2}$, see figure \ref{Newtoninterpolation},
by corollary \ref{Corollaryoflemma}, we have the following identity
$$\begin{array}{rl}
x^n&=e(a_0,b)\\
& \quad +e(a_1,b)\cdot(x-x_1)\\
& \qquad+e(a_2,b)\cdot(x-x_1)(x-x_2)+\cdots \\
& \qquad \quad +e(a_n,b)\cdot(x-x_1)\cdots(x-x_n),
\end{array}$$
with $a_i=(i,1)$, $N=(n,n)$.

One can show by induction that for each point $(i,j)$ with $0\leq i,j\leq n$, $e((i,j),b)$
corresponds to a divided differences of $x^n$.
More exactly, consider the inverse diagram, and stare at the following diagram.
$$\xymatrix{
\vdots\ar[d] & \vdots\ar[d] & \\
x_2^n\ar[d]& \diamondsuit\ar[l]^{(x_2-x_3)}\ar[d]&\ar[l]\ar[d]\cdots\\
x_1^n& \heartsuit\ar[l]^{(x_1-x_2)}&\ar[l]^{(x_1-x_3)}\square& \ar[l]\cdots}$$
Since $\heartsuit\cdot (x_1-x_2)+x_2^n=x_1^n$, so $\heartsuit=\frac{x_1^n-x_2^n}{x_1-x_2}$.
The same reason, $\diamondsuit=\frac{x_2^2-x_3^n}{x_2-x_3}$. Therefore $\square = \frac{\heartsuit-\diamondsuit}{x_1-x_3}$, and so on.

\vfill
\bibliographystyle{plain}
\bibliography{bibfile}
\end{document}